\documentclass[11pt]{article}
\usepackage{mathrsfs}
\usepackage{amssymb}
\usepackage{amsmath}
\usepackage[all]{xy}

\def\Ker{\mathop{\rm Ker}\nolimits}
\def\Hom{\mathop{\rm Hom}\nolimits}
\def\Coker{\mathop{\rm Coker}\nolimits}
\def\Mod{\mathop{\rm Mod}\nolimits}

\def\rfd{\mathop{\rm r.fd}\nolimits}
\def\End{\mathop{\rm End}\nolimits}

\def\max{\mathop{\rm max}\nolimits}

\def\rfd{\mathop{\rm r.fd}\nolimits}

\title{\Large \bf The Auslander-Type Condition of Triangular Matrix Rings
}
\author{Chonghui Huang$^{1,2}$, Zhaoyong Huang$^1$\thanks{{\it E-mail address}: huangzy@nju.edu.cn}
\\\\{\small \it $^{1}$Department of Mathematics, Nanjing University,
Nanjing 210093, China;} \\ {\small \it $^{2}$College of Mathematics
and Physics, University of South China, Hengyang 421001, China}}

\date{}
\begin{document}
\baselineskip=18pt \maketitle

\begin{abstract}
Let $R$ be a left and right Noetherian ring and $n,k$ any
non-negative integers. $R$ is said to satisfy the Auslander-type
condition $G_n(k)$ if the right flat dimension of the $(i+1)$-st
term in a minimal injective resolution of $R_R$ is at most $i+k$ for
any $0 \leq i \leq n-1$. In this paper, we prove that $R$ is
$G_n(k)$ if and only if so is a lower triangular matrix ring of any
degree $t$ over $R$.
\end{abstract}

\vspace{0.5cm}

\centerline{\bf 1. Introduction}

\vspace{0.2cm}

Let $R$ be a ring and $M$ a right $R$-module. We use $$0\to M\to
I^0(M) \to I^1(M) \to \cdots \to I^i(M) \to \cdots$$ to denote a
minimal injective resolution of $M_R$. For a positive integer $n$,
recall from [FGR] that a left and right Noetherian ring $R$ is
called an {\it $n$-Gorenstein ring} if the right flat dimension of
$I^i(R)$ is at most $i$ for any $0 \leq i \leq n-1$, and $R$ is said
to satisfy the {\it Auslander condition} if $R$ is $n$-Gorenstein
for all $n$. The notion of the Auslander condition may be regarded
as a non-commutative version of commutative Gorenstein rings. A
remarkable property of $n$-Gorenstein rings (and hence rings
satisfying the Auslander condition) is its left-right symmetric,
which was proved by Auslander (see [FGR, Theorem 3.7]). Motivated by
the philosophy of Auslander, Huang and Iyama introduced in [HuIy]
the notion of the Auslander-type condition as follows. For any
$n,k\geq 0$, a left and right Noetherian ring $R$ is said to be
$G_n(k)$ if the right flat dimension of $I^i(R)$ is at most $i+k$
for any $0 \leq i \leq n-1$. It is trivial that $R$ is an
$n$-Gorenstein ring if and only if $R$ is $G_n(0)$. In general case,
the Auslander-type condition $G_n(k)$ does not possess the
left-right symmetry (see [HuIy]). Note that the Auslander-type
condition plays a crucial role in representation theory of algebras
and homological algebra (e.g. [AuR1,2], [B], [EHIS], [FGR], [FI],
[HN], [Hu1,2], [HuIy], [I], [IM], [IS], [IW], [Iy1,2,3,4], [M],
[W]).

It was proved by Iwanaga and Wakamatsu in [IW, Theorem 8] that a
left and right Artinian ring $R$ is an $n$-Gorenstein ring if and
only if so is a lower triangular matrix ring of any degree $t$ over
$R$. Observe that this is a generalization of [FGR, Theorem 3.10]
where the case $k=2$ was established. In this paper, we will
generalize the Iwanaga and Wakamatsu's result mentioned above, and
prove the following result.

\vspace{0.2cm}

{\bf Theorem} {\it Let $R$ be a left and right Noetherian ring and
$n,k\geq 0$. Then $R$ is $G_n(k)$ if and only if so is a lower
triangular matrix ring $T_t(R)$ of any degree $t$ over $R$.}

\vspace{0.2cm}

In Section 2, we recall some notions and notations and give some
preliminary results about triangular matrix rings. Then in Section
3, we give the proof of the above theorem, by establishing the
relation between the flat dimensions of the corresponding terms in
the minimal injective resolutions of $R_R$ and $T_t(R)_{T_t(R)}$. In
[Iy1], Iyama introduced the notions of the $(l,n)^{op}$-condition
(which has a close relation with the Auslander-type condtion) and
the dominant number. In Section 3, we also prove the following
results. Let $R$ be a left and right Noetherian ring and $l,n \geq
0$, $t\geq 1$. If $R$ satisfies the $(l, n)^{op}$-condition, then
$T_t(R)$ satisfies the $(l+1, n)^{op}$-condition. Conversely, if
$T_t(R)$ satisfies the $(l, n)^{op}$-condition, then so does $R$. In
addition,  If $n$ is a dominant number of $R$, then $n+1$ is a
dominant number of $T_t(R)$.

\vspace{0.5cm}

\centerline{\bf 2. Preliminaries}

\vspace{0.2cm}

In this section, we give some notions and notations and collect some
elementary facts which are useful for the rest of this paper.

Throughout this paper, $R$ and $S$ are rings and $_SM_R$ is a left
$S$ right $R$-bimodule. We denote by $\Lambda$=$\left(
\begin{array}{cc}
R & 0 \\
M & S \\
\end{array}
\right)$ the triangular matrix ring, and denote by $^*(-)$ the
functor $\Hom_R(M, -)$. For the ring $R$, we use $\Mod R$ to denote
the category of right $R$-modules.

By [G], $\Mod \Lambda$ is equivalent to a category $\mathscr{T}$ of
triples $(X,Y)_f$, where $X\in \Mod R$ and $Y\in \Mod S$ and $f:\
Y\otimes_S M_R \to X_R$ is a homomorphism in $\Mod R$ (which is
called the {\it associated homomorphism}). The right
$\Lambda$-module corresponding to the triple $(X,Y)_f$ is the
additive group $X\oplus Y$ with the right $\Lambda$-action given by
$$(x,y)\left(
\begin{array}{cc}
r & 0 \\
m & s \\
\end{array}
\right)=(xr+f(y\otimes m), ys)$$ for any $x\in X, y\in Y, r\in R,
s\in S$ and $m\in M$.

Another description of a right $\Lambda$-module $X\oplus Y$ is a
triple $_{\varphi}(X,Y)$, where $\varphi: \ \ Y_S\to \Hom_R(_SM_R,
X_R)_S$ is a homomorphism in $\Mod S$ (which is also called the {\it
associated homomorphism}). The operation of the additive group
$X\oplus Y$ with the right $\Lambda$-action given by
$$(x,y)\left(
\begin{array}{cc}
r & 0 \\
m & a \\
\end{array}
\right)=(xr+\varphi(y)(m), ys)$$ for any $x\in X, y\in Y, r\in R,
s\in S$ and $m\in M$.




In particular, we have the following isomorphism:
$$\Hom_S(Y_S, \Hom_R(_SM_R, X_R)_S)\cong \Hom_R(Y\otimes_SM_R,
X_R).$$ So it is convenient for us to adopt either of these two
descriptions of $X\oplus Y$ in the following argument.

If $(U,V)_g$ and $(X,Y)_f$ are in $\mathscr{T}$, then the
homomorphism from $(U,V)_g$ to $(X,Y)_f$ are pairs $(h_1, h_2)$,
where $h_1:\ U\to X$ is a homomorphism in $\Mod R$ and $h_2: \ V\to
Y$ is a homomorphism in $\Mod S$ satisfying the condition $h_1
g=f(h_2\otimes 1_M)$. If $M\in \Mod S$ is flat, then it is not
difficult to verify that $(h_1, h_2)$ is monic (resp. epic) if and
only if so are both of $h_1$ and $h_2$.

\vspace{0.2cm}

{\bf Lemma 2.1} ([FGR, Proposition 1.14]) {\it Let $X$, $Y$ and $f$
be as above. Then $(X, Y)_f\in \Mod \Lambda$ is flat if and only if
the following conditions are satisfied.

(1) $Y \in \Mod S$ is flat.

(2) $\Coker f \in \Mod R$ is flat.

(3) $f$ is a monomorphism.}

\vspace{0.2cm}

{\bf Lemma 2.2} ([Y, Corollary 6]) {\it $\Lambda$ is a left (resp.
right) Noetherian ring if and only if both $R$ and $S$ are left
(resp. right) Noetherian rings and $_SM$ (resp. $M_R$) is finitely
generated}.

\vspace{0.2cm}

For the ring $R$ and any positive integer $t$, we use $T_t(R)$ to
denote the triangular matrix ring $\left(
\begin{array}{cccc}
R & \  &\   &\ \\
R & R& \ & \ \\
\vdots & \vdots &\ddots & \ \\
R & R & \cdots & R \\
\end{array}
\right)$ of degree $t$.

\vspace{0.2cm}

{\bf Lemma 2.3} {\it For any $t \geq 2$, $T_t(R)$ is a triangular
matrix ring of the form: $$T_t(R)=\left(
\begin{array}{cc}
T_{t-1}(R) & 0 \\
_RR^{(t-1)}_{T_{t-1}(R)} & R \\
\end{array}
\right).$$ In particular, $R^{(t-1)}_{T_{t-1}(R)}$ is faithful and
finitely generated projective and $\End
_{T_{t-1}(R)}(R^{(t-1)})\cong R$.}

\vspace{0.2cm}

{\bf Proof.} We can regard $R^{(t-1)}$ as a right
$T_{t-1}(R)$-module in a natural way. Let $e=\left(
\begin{array}{cccc}
0 & 0 &\cdots   &0 \\
0 & 0& \cdots & 0 \\
\vdots & \vdots &\  & \vdots \\
0& 0 & \cdots & 1 \\
\end{array}
\right)$ be a matrix in $T_{t-1}(R)$ such that the $(t-1,
t-1)$-component is 1 and 0 elsewhere. Then $e$ is an idempotent and
$R^{(t-1)}_{T_{t-1}(R)}\cong eT_{t-1}(R)_{T_{t-1}(R)}$, which
implies that $R^{(t-1)}_{T_{t-1}(R)}$ is faithful and finitely
generated projective and $\End _{T_{t-1}(R)}(R^{(t-1)})\cong
eT_{t-1}(R)e\cong R$. \hfill $\square$

\vspace{0.2cm}

{\bf Proposition 2.4} {\it If $R$ is a left (resp. right) Noetherian
ring, then so is $T_t(R)$ for any $t\geq 1$.}

\vspace{0.2cm}

{\bf Proof}. We proceed by induction on $t$. The case for $t=1$ is
trivial, and the case for $t=2$ follows from Lemma 2.2. Now suppose
$t\geq 3$. By Lemma 2.3, $T_t(R)=\left(
\begin{array}{cc}
T_{t-1}(R) & 0 \\
_RR^{(t-1)}_{T_{t-1}(R)} & R \\
\end{array}
\right)$ with both $_RR^{(t-1)}$ and $R^{(t-1)}_{T_{t-1}(R)}$
finitely generated. Then by the induction hypothesis and Lemma 2.2,
we get the assertion. \hfill $\square$

\vspace{0.2cm}

{\bf Definition 2.5} ([EJ]) Assume that $\mathscr{F}$ is a subclass
of $\Mod R$, $X\in\mathscr{F}$ and $Y\in \Mod R$. The homomorphism
$f: X\to Y$ is said to be an $\mathscr{F}$-{\it precover} of $Y$ if
$\Hom _{R}(X', X)\to \Hom _{\Lambda}(X', Y)\to 0$ is exact for any
$X'\in\mathscr{F}$. An $\mathscr{F}$-precover $f: X\to Y$ is said to
be an $\mathscr{F}$-{\it cover} of $Y$ if an endomorphism $g:X\to X$
is an automorphism whenever $f=fg$. If $\mathscr{F}$ is the subclass
of $\Mod R$ consisting of all flat right $R$-modules, then an
$\mathscr{F}$-cover is called a {\it flat cover}.

\vspace{0.2cm}

Bican, El Bashir and Enochs proved in [BEE, Theorem 3] that every
module in $\Mod R$ has a flat cover. For a module $N \in \Mod R$, we
call the following exact sequence:
$$\cdots \to F_i(N)\buildrel {{\pi}_i}(N) \over \to \cdots \to F_1(N)
\buildrel {{\pi}_1}(N) \over\to F_0(N) \buildrel {{\pi}_0}(N) \over\to N \to 0$$
a {\it minimal flat resolution} of $N_R$, where ${{\pi}_0}(N):\
F_0(N) \to N$ is a flat cover of $N$ and ${{\pi}_i}(N):\ F_i(N) \to
\Ker {{\pi}_{i-1}}(N)$ is a flat cover of $\Ker {{\pi}_{i-1}}(N)$
for any $i \geq 1$. We denote the right flat dimension of $N$ by
$\rfd _R(N)$. It is easy to verify that $\rfd _R(N)\leq n$ if and
only if $F_{n+1}(N)=0$.

\vspace{0.5cm}

\centerline{\bf 3. Main Results}

\vspace{0.2cm}

In this section, we give the proof of the main result mentioned in
the Introduction, by establishing the relation between the flat
dimensions of the corresponding terms in the minimal injective
resolutions of $R_R$ and $T_t(R)_{T_t(R)}$.

\vspace{0.2cm}

{\bf Lemma 3.1} {\it Let $X$, $Y$ and $f$ be as in Section 2. If $M$
is left $S$-flat, then $(F_0(Y)\otimes_S M, F_0(Y))_1\oplus
(F_0(X),0)_0 \buildrel {\psi _0} \over \longrightarrow (X,Y)_{f}\to
0$ is an exact sequence in $\Mod \Lambda$ with $(F_0(Y)\otimes_S M,
F_0(Y))_1\oplus (F_0(X),0)_0$ flat and $\psi _0=((f
(\pi_0(Y)\otimes_S 1_M), \pi_0(X)), \pi_0(Y))$. Moreover, if $f_1$
is the associated homomorphism of $\Ker \psi _0$, then $ 0\to \Ker f
\to \Coker f_1 \to F_0(X) \to \Coker f \to 0$ is an exact sequence
in $\Mod R$.}

\vspace{0.2cm}

{\bf Proof}. Since $M$ is left $S$-flat, we have the following
commutative diagram with exact rows:
$$\xymatrix{
0 \ar[r]& \Ker \pi_0(Y)\otimes_S M \ar[d]^{f_1}
\ar[r]&F_0(Y)\otimes_S M \ar[d]^{\binom {1_{F_0(Y)\otimes_S M}} 0}
\ar[r]^{\pi_0(Y)\otimes_S 1_M}& Y\otimes_S M\ar[d]^{f} \ar[r]& 0  \\
0 \ar[r] & \Ker h  \ar[r] & (F_0(Y)\otimes_S M)\oplus F_0(X)
 \ar[r]^{\ \ \ \ \ \ \ \ \ \ h} & X
 \ar[r] & 0  }$$ where $h=(f(\pi_0(Y)\otimes_S 1_M),
\pi_0(X))$ and $f_1$ is established by diagram-chasing. Then by
Lemma 2.1, ${((F_0(Y)\otimes_S M), F_0(Y))_1\oplus(F_0(X),0)_0}\in
\Mod \Lambda$ is flat. The last assertion follows from the snake
lemma. \hfill$\square$

%

\vspace{0.2cm}

From now on, assume that $M_R$ is finitely generated, faithful and
projective with $S=\End_R(M)$ and $_SM$ is finitely generated
projective. Then by [IW, Corollary 3],
$$I^0(\Lambda)={_1(I^0(R), {^*(I^0(R))})}\oplus
 {_1(I^0(M), {^*(I^0(M))})}$$  and  $$I^i(\Lambda)={_1(I^i(R),
 {^*(I^i(R))})}\oplus {_1(I^i(M), {^*(I^i(M))})}\oplus {_0(0, I^{i-1}({^*R}))}\ \ \ \ \
(i\geq 1)$$ give a minimal injective resolution of $\Lambda
_{\Lambda}$. In the following, we will construct a flat resolution
of $I^i(\Lambda_{\Lambda})$ for any $i \geq 0$, and then consider
the Auslander-type condition of the triangular matrix ring $T_t(R)$.

\vspace{0.2cm}

{\bf Proposition 3.2} {\it (1) Let $I_R$ be injective and $\xi_{I}:
{^{*}I}\otimes_S M \to I$ defined by $\xi_{I}(\alpha\otimes
x)=\alpha(x)$ for any $\alpha\in {^{*}I}$ and $x\in M$ be the
natural homomorphism. Then
\begin{eqnarray*}
\nonumber F_0&=&{(F_0({^*I})\otimes_S M, F_0({^*I}))}_1, \\
\nonumber \Ker \psi _0&=&{(\Ker \xi_I(\pi_0({^*I})\otimes_S 1_M), \Ker \pi_0({^*I}))}_{f_1}, \\
\nonumber F_i& =& {(F_i({^*I})\otimes_S M, F_i({^*I}))}_1
\oplus (F_0(\Ker f_{i-1}(\pi_{i-1}({^*I})\otimes_S 1_M)), 0)_0, \\
\nonumber \Ker \psi _i&=&{(\Ker f_i(\pi_i({^*I})\otimes_S 1_M), \Ker
\pi_i({^*I}))}_{f_{i+1}}\ \ \ \ \ (i\geq 1)
\end{eqnarray*} give a flat resolution of an injective right $\Lambda$-module $_1(I_R,
{^*I_R})$:
$$\cdots \to F_n \buildrel {\psi _n} \over \longrightarrow F_{n-1} \cdots \to F_1\buildrel {\psi _1}
\over \longrightarrow F_0 \buildrel {\psi _0} \over \longrightarrow
\ _1(I_R, {^*I_R})\to 0,$$ where $f_i $ is established by
diagram-chasing as in Lemma 3.1. In particular, $\rfd_{\Lambda}
{_1(I_R, {^*I_R})}\leq k$ if and only if $\rfd_R \Ker \xi_I \leq
k-1$ and $\rfd_S {^*I_R}\leq k$.

(2) If $\Hom _R({_SM_R}, R)$ is finitely generated right
$S$-projective and $E\in \Mod S$, then
\begin{eqnarray*}
\nonumber F_0&=&{(F_0(E)\otimes_S M, F_0({E}))}_1 \ and\\
\nonumber  F_i&=& {(F_i({E})\otimes_S M, F_i({E}))}_1\oplus
(F_{i-1}(E)\otimes_S M, 0)_0\ \ \ \ \ (i\geq 1)
\end{eqnarray*} give a flat resolution of $_0(0, E_S)$ in $\Mod \Lambda$. In particular, $\rfd_{\Lambda}
{_0(0, E_S)} \leq k$ if and only if $\rfd_S E\leq k-1$.}


\vspace{0.2cm}

{\bf Proof}. (1) We proceed by induction on $i$. Since $M_R$ is
finitely generated, faithful and projective with $S=\End_R(M)$, by
[AF, Proposition 20.11], it is not difficult to verify that
$\xi_{I}$ is epic. Thus we have the following exact sequence:
$$0\to \Ker \xi_{I} \to {^*I}\otimes_S M
\buildrel {\xi_{I}} \over \to I\to 0.$$ Then $F_0$ and $\Ker \psi
_0$ are given by the following commutative diagram with exact rows:
$$\xymatrix{ 0 \ar[r]& \Ker
\pi_0({^*I})\otimes_S M \ar[d]^{f_1} \ar[r]&F_0({^*I})\otimes_S M
\ar@{=}[d] \ar[r]^{\ \ \ \ \pi_0({^*I})\otimes_S 1_M}& \ \ \ \ \ \
{^*I}\otimes_S M \ \ \ \  \ \ar[d]^{\xi_I} \ar[r]& 0  \\
0 \ar[r] & \Ker \xi_{I}(\pi_0({^*I})\otimes_S 1_M) \ar[r] &
F_0({^*I})\otimes_S M \ar[r]^{\ \ \ \  \ \ \ \ \ \ \
\xi_{I}(\pi_0({^*I})\otimes_S 1_M)} &\ \ \  I \ar[r] & 0  }$$ where
$f_1$ is established by diagram-chasing. By the snake lemma, we have
an exact sequence $0\to \Ker \pi_0({^*I})\otimes_S M \buildrel
{f_{1}} \over \longrightarrow \Ker \xi_{I}(\pi_0({^*I})\otimes_S M)
\to \Ker \xi_I\to 0$. Then by using Lemma 3.1 iteratively and the
induction hypothesis, we get the homomorphism $f_{i+1}$ (which is
monic), $F_i = {(F_i({^*I})\otimes_S M, F_i({^*I}))}_1 \oplus
(F_0(\Ker f_{i-1}(\pi_{i-1}({^*I})\otimes_S M)), 0)_0$ and $\Ker
\psi _i={(\Ker f_i(\pi_i({^*I})\otimes_S 1_M), \Ker
\pi_i({^*I}))}_{f_{i+1}}$; in particular, we get an exact sequence
$$0\to \Coker f_{i+1} \to F_0(\Ker f_i(\pi_i({^*I})\otimes_S 1_M))\to
\Coker f_i \to 0.$$ So $\Coker f_i$ is right $R$-flat if and only if
$\rfd_R \Ker \xi_I\leq i-1$. In addition, by Lemma 2.1, $_1(I_R,
{^*I_R})$ is right $\Lambda$-flat if and only if ${^*I_R}$ is right
$S$-flat and $\Ker \xi_I=0$. So for any $k \geq 0$, we have that
$\rfd_{\Lambda} {_1(I_R, {^*I_R})}\leq k$ if and only if $\Ker \psi
_{k-1}$ is right $\Lambda$-flat, if and only if $\Coker f_k$ is
right $R$-flat and $\Ker \pi _{k-1}(^{*}I_R)$ is right $S$-flat, if
and only if $\rfd_R \Ker \xi_I \leq k-1$ and $\rfd_S {^*I_R}\leq k$.

(2) By Lemma 3.1, we have that $F_0={(F_0(E)\otimes_S M,
F_0({E}))}_1$ and there exists an exact sequence: $$0\to \Ker
\pi_0(E)\otimes_S M \buildrel {f_1} \over \longrightarrow
F_0(E)\otimes_S M \to E\otimes_S M\to 0.$$ By using an argument
similar to that in (1), we have that $F_i= {(F_i({E})\otimes_S M,
F_i({E}))}_1\oplus (F_{i-1}(E)\otimes_S M, 0)_0$ for any $i\geq 1$,
and that $\Coker f_i$ is right $R$-flat if and only if $\rfd_R
E\otimes_S M \leq i-1$. In addition, ${_0(0, E_S)}$ is right
$\Lambda$-flat if and only if $E_S$ is right $S$-flat. Thus for any
$k \geq 0$, $\rfd_{\Lambda} {_0(0, E_S)} \leq k$ if and only if
$\rfd _SE\leq k$ and $\Coker f_k$ is right $R$-flat, if and only if
$\rfd_S E\leq k$ and $\rfd_R E\otimes_S M \leq k-1$. So, to get the
second assertion, it suffices to prove that $\rfd_R E\otimes_S M\leq
k-1$ if and only if $\rfd_S E\leq k-1$.

Let $$\cdots \to F_i\buildrel {f_i} \over \longrightarrow F_{i-1}
\to \cdots \to F_1 \buildrel {f_1} \over \longrightarrow
F_0\buildrel {f_0} \over \longrightarrow E \to 0$$ be a flat
resolution of $E_S$ in $\Mod S$. Since $_SM$ is projective, $$\cdots
\to F_i\otimes_S M\buildrel {f_i\otimes_S 1_M} \over \longrightarrow
F_{i-1}\otimes_S M \to \cdots \to F_1\otimes_S M \buildrel
{f_1\otimes_S 1_M} \over \longrightarrow F_0\otimes_S M\buildrel
{f_0\otimes_S 1_M} \over \longrightarrow E\otimes_S M \to 0$$ is a
flat resolution of $E\otimes _SM$ in $\Mod R$. So, if $\rfd_S E\leq
k-1$, then $\rfd_R E\otimes_S M\leq k-1$.

Conversely, assume that $\rfd_R E\otimes_S M\leq k-1$. Then $\Coker
f_k\otimes _S M$ is right $R$-flat. So $\Coker f_k\otimes _S M$ is a
direct limit of a direct system of finitely generated projective
right $R$-modules $\{ Q_i\}_{i\in I}$, that is, $\Coker f_k\otimes
_S M=\underset{i\in I}{\underset{\rightarrow}{\rm lim}}Q_i$, where
$I$ is a direct index set. Because $M_R$ is finitely generated
projective and $S=\End _R(M)$, by [AF, Proposition 20.10] and [GT,
Lemma 1.2.5], we have that $$\Coker f_k\cong \Hom_R({_SM_R}, \Coker
f_k\otimes_S M)\cong
 \Hom_R({_SM_R}, \underset{i\in I}{\underset{\rightarrow}{\rm
lim}}Q_i)\cong\underset{i\in I}{\underset{\rightarrow}{\rm
lim}}\Hom_R({_SM_R}, Q_i).$$ By assumption, $\Hom_R({_SM_R}, R)$ is
finitely generated right $S$-projective, so $\Coker f_k$ is right
$S$-flat and $\rfd_S E\leq k-1$. \hfill $\square$

\vspace{0.2cm}

{\bf Proposition 3.3} {\it If $\Hom _R({_SM_R}, R)$ is finitely
generated right $S$-projective, then for any $k, i\geq 0$,
$\rfd_\Lambda I^i(\Lambda)\leq k$ if and only if the following
conditions are satisfied.

(1) $\rfd_R \Ker\xi_{I^i(R)} \leq k-1$.

(2) $\rfd_S {^*(I^i(R))}\leq k$.

(3) $\rfd_S I^{i-1}({^*R})\leq k-1$ (Here we set $I^{-1}({^*R})=0$
and $\rfd _S 0 = -1 $).}

\vspace{0.2cm}

{\bf Proof}. Since $M_R$ is finitely generated projective,
$\rfd_\Lambda {_1(I^i(R), {^*(I^i(R))})}\leq k$ yields that
$\rfd_\Lambda {_1(I^i(M), {^*(I^i(M))})}\leq k$. Note that
$\rfd_\Lambda I^0(\Lambda)\leq k$ if and only if $\rfd_\Lambda
{_1(I^0(R), {^*(I^0(R))})}$\linebreak $\leq k$. So, by Proposition
3.2(1), $\rfd_\Lambda I^0(\Lambda)\leq k$ if and only if $\rfd_R
\Ker\xi_{I^0(R)} \leq k-1$ and $\rfd_S {^*(I^0(R))}\leq k$. The case
for $i=0$ follows. Now suppose $i\geq 1$. Note that $\rfd_\Lambda
I^i(\Lambda)\leq k$ if and only if $\rfd_\Lambda {_1(I^i(R),
{^*(I^i(R))})}\leq k$ and $\rfd_\Lambda {_0(0, I^{i-1}({^*R}))}\leq
k$. So, by Proposition 3.2, we have that $\rfd_\Lambda
I^i(\Lambda)\leq k$ if and only if $\rfd_R \Ker\xi_{I^i(R)} \leq
k-1$, $\rfd_S {^*(I^i(R))}\leq k$ and $\rfd_S I^{i-1}({^*R})\leq
k-1$. \hfill$\square$

\vspace{0.2cm}

{\bf Proposition 3.4} {\it If $\Hom _R({_SM_R}, R)$ is finitely
generated right $S$-projective, then for any $i\geq 0$,
$\rfd_\Lambda I^i(\Lambda)\leq k$ yields that $\rfd_R I^i(R)\leq k$
and $\rfd_S I^{i-1}(S)\leq k-1$}.

\vspace{0.2cm}

{\bf Proof}. For any $i\geq 0$, by Proposition 3.3(1) and (2), we
have that $\rfd_R \Ker \xi_{I^i(R)}\leq k-1$ and $\rfd_R
{^*(I^i(R))}\otimes_S M\leq k$. In addition, we have the following
exact sequence:
$$0\to \Ker \xi_{I^i(R)} \to {^*(I^i(R))}\otimes_S M
\buildrel {\xi_{I^i(R)}} \over \longrightarrow I^i(R)\to 0,$$ which
implies $\rfd_R I^i(R)\leq k$.

On the other hand, since $M_R$ is finitely generated projective, the
condition (3) in Proposition 3.3 is also satisfied when $^{*}R$ is
replaced by $^{*}M$. It follows that $\rfd_S I^{i-1}(S)\leq k-1$.
\hfill$\square$

\vspace{0.2cm}

{\bf Proposition 3.5} {\it For any $i\geq 0$, $\rfd_\Lambda
I^i(\Lambda)\leq k$ if and only if $\rfd_\Gamma I^i(\Gamma)\leq k$,
where $\Gamma =\left(
\begin{array}{ccc}
R & 0 &0 \\
M & S & 0 \\
M & S & S \\
\end{array}
\right)$.}

\vspace{0.2cm}

{\bf Proof}. Let $e=\left(
\begin{array}{cc}
0 & 0 \\
0 & 1 \\
\end{array}
\right)\in \Lambda$. Then we have $$\Gamma=\left(
\begin{array}{cc}
\Lambda & 0 \\
e\Lambda & e\Lambda e \\
\end{array}
\right),\ \Lambda e=e\Lambda e\cong S.$$

Since $\Lambda$ can be embedded in
$$\End_S(e\Lambda)\cong \left(
\begin{array}{cc}
\End_S(M) & \Hom_S(M, S) \\
M & S \\
\end{array}
\right),$$ $e\Lambda$ is a faithful right $\Lambda$-module. It is
trivial that $e\Lambda$ is finitely generated projective as a right
$\Lambda$-module and a left $S$-module. Notice that $S\cong
\End_\Lambda(e\Lambda)$ and $\Hom_\Lambda(e\Lambda, \Lambda)_S \cong
\Lambda e_S\cong S_S$, so by Proposition 3.3, we get that
$\rfd_\Gamma I^i(\Gamma)\leq k$ if and only if the following
conditions are satisfied.

(1) $\rfd_\Lambda \Ker \eta _{I^i(\Lambda)} \leq k-1$, where
$\eta_{E}: \Hom_\Lambda(e\Lambda, E)\otimes_{e\Lambda e}e\Lambda \to
E$ defined by $\eta _{E}(\alpha\otimes x)=\alpha(x)$ for any
$\alpha\in \Hom_\Lambda(e\Lambda, E)$ and $x\in e\Lambda$ is the
natural homomorphism for an injective right $R$-module $E$.

(2) $\rfd_S \Hom_\Lambda(e\Lambda, I^i(\Lambda))_S\leq k$.

(3) $\rfd_S I^{i-1}(\Hom_\Lambda(e\Lambda, \Lambda)_S)\leq k-1$.

Assume that $\rfd_\Gamma I^i(\Gamma)\leq k$. By Proposition 3.4,
$\rfd_\Lambda I^i(\Lambda)\leq k$. It remains to show that
$\rfd_\Lambda I^i(\Lambda)\leq k$ implies $\rfd_\Gamma
I^i(\Gamma)\leq k$ for any $i \geq 0$. To do this, it suffices to
show the conditions $(1)-(3)$ above are satisfied.

Note that $\Hom_\Lambda(e\Lambda, \Lambda)_S \cong S_S$. If $L\in
\Mod \Lambda$ is flat, then $L$ is a direct limit of a direct system
of finitely generated projective right $\Lambda$-modules $\{
P_i\}_{i\in I}$, that is, $L=\underset{i\in
I}{\underset{\rightarrow}{\rm lim}}P_i$, where $I$ is a direct index
set. So $\Hom _{\Lambda}(e\Lambda, L) \cong \Hom_{\Lambda}(e\Lambda,
\underset{i\in I}{\underset{\rightarrow}{\rm lim}}P_i)\cong
\underset{i\in I}{\underset{\rightarrow}{\rm
lim}}\Hom_{\Lambda}(e\Lambda, P_i)$ is right $S$-flat. Then it is
not difficult to verify that $\rfd_\Lambda I^i(\Lambda)\leq k$
yields $\rfd_S \Hom_\Lambda(e\Lambda, I^i(\Lambda))_S\leq k$. Thus
the condition $(2)$ is satisfied.

Notice that $\Hom_\Lambda(e\Lambda, \Lambda)_S \cong S_S$, so by
Proposition 3.4, the condition $(3)$ is satisfied.

By [IW, Corollary 3], $$I^0(\Lambda)={_1(I^0(A),
{^*(I^0(R))})}\oplus {_1(I^0(M), {^*(I^0(M))})},$$ and
$$I^i(\Lambda)={_1(I^i(R), {^*(I^i(R))})}\oplus {_1(I^i(M),
{^*(I^i(M))})}\oplus {_0(0, I^{i-1}({^*R}))}$$ give a minimal
injective resolution of $\Lambda_\Lambda$. So, to verify the
condition (1), it suffices to show that $\rfd_\Lambda \Ker \eta _{E}
\leq k-1$ for any injective right $\Lambda$-module $E$, where $E$ is
of the form:

(a) ${_0(0, I_S)}$ with $I_S$ injective, or

(b) ${_1(I_R, {^*(I_R)})}$ with $I_R$ injective.

If $E$ is of the form (a), then $$\Hom_\Lambda(e\Lambda,
E)\otimes_{e\Lambda e}e\Lambda\cong Ee\otimes_S e\Lambda\cong
E\otimes_S (M, S)_1\cong E,$$ which implies $\Ker \eta_E=0$.

If $E$ is of the form (b), then $E={_1(I_R, {^*(I_R)})}$ and
$$\Ker \eta_E\cong {_0(\Ker \xi_{I_R}, 0)}.$$ On the other hand, by
Lemma 3.1, we have $\rfd_\Lambda {_0(X_R, 0)} =\rfd_R X$. So
$\rfd_\Lambda \Ker \eta_{I^i(\Lambda)} \leq k-1$ if and only if
$\rfd_R \Ker \xi_{I^i(R)} \leq k-1$. Because $\rfd_\Lambda
I^i(\Lambda)\leq k$, $\rfd_\Lambda \Ker \eta _{E} \leq k-1$. \hfill
$\square$

\vspace{0.2cm}

We are now in a position to state the main result in this section.

\vspace{0.2cm}

{\bf Theorem 3.6} {\it If $R$ is a left and right Noetherian ring
and $t$ is a positive integer, then $T_t(R)$ is a left and right
Noetherian ring, and $\rfd_{T_t(R)} I^i(T_t(R))=\max\{\rfd_R
I^i(R),$ $\rfd_R I^{i-1}(R)+1\}$ for any $i\geq 0$.}

\vspace{0.2cm}

{\bf Proof}. The first assertion follows from Proposition 2.4. We
will prove the second assertion by induction on $t$. The case for
$t=1$ is trivial, and the case for $t=2$ follows from [FGR, Theorem
3.10].

Now suppose $t\geq 3$. By Lemma 2.3, $T_t(R)=\left(
\begin{array}{cc}
T_{t-1}(R) & 0 \\
_RR^{(t-1)}_{T_{t-1}(R)} & R \\
\end{array}
\right)$ with $R^{(t-1)}_{T_{t-1}(R)}$ faithful and finitely
generated projective and $\End _{T_{t-1}(R)}(R^{(t-1)})\cong R$.
Then by Proposition 3.5, $\rfd_{T_t(R)} I^i(T_t(R))\leq k$ if and
only if $\rfd_{T_{t-1}(R)} I^i(T_{t-1}(R))\leq k$. So we have that
$\rfd_{T_t(R)} I^i(T_t(R))$=$\rfd_{T_2(R)} I^i(T_2(R))= \max\{\rfd_R
I^i(R), \rfd_R I^{i-1}(R)+1\}$ by the induction hypothesis. \hfill
$\square$

\vspace{0.2cm}

As an immediate consequence of Theorem 3.6, we get the main theorem
mentioned in the Introduction.

\vspace{0.2cm}

{\bf Theorem 3.7} {\it If $R$ is a left and right Noetherian ring
and $n, k\geq 0$, $t \geq 1$, then $R$ is $G_n(k)$ if and only if so
is $T_t(R)$.}

\vspace{0.2cm}

{\bf Proof.} By Theorem 3.6, we have that $\rfd_R I^i(R)\leq
\rfd_{T_t(R)} I^i(T_t(R))$ for any $i\geq 0$, so the sufficiency is
trivial. Conversely, if $R$ is $G_n(k)$, then by Theorem 3.6,
$\rfd_{T_t(R)} I^i(T_t(R))=\max\{\rfd_R I^i(R), \rfd_R
I^{i-1}(R)+1\}\leq i+k$ for any $0 \leq i\leq n-1$ and $T_t(R)$ is
$G_n(k)$.  \hfill $\square$

\vspace{0.2cm}

We recall some notions introduced by Iyama in [Iy1]. Let $R$ be a
left and right Noetherian ring and $l, n\geq 0$. $R$ is said to
satisfy the {\it $(l,n)^{op}$-condition} if $\rfd I_i(R) \leq l-1$
for any $0 \leq i \leq n-1$. It is easy to see that $R$ is $G_n(k)$
if and only if $R$ satisfies the $(k+i,i)^{op}$-condition for any $1
\leq i \leq n$. In addition, if $\rfd_R I^i(R)< \rfd_R I^n(R)$ for
any $0 \leq i \leq n-1$, then $n$ is called a {\it dominant number}
of $R_R$. As another application of Theorem 3.6, we get the
following

\vspace{0.2cm}

{\bf Corollary 3.8} {\it If $R$ is a left and right Noetherian ring,
then for any $l, n\geq 0$, $t \geq 1$, we have the following

(1) If $R$ satisfies the $(l, n)^{op}$-condition, then $T_t(R)$
satisfies the $(l+1, n)^{op}$-condition. Conversely, if $T_t(R)$
satisfies the $(l, n)^{op}$-condition, then so does $R$.

(2) If $n$ is a dominant number of $R$, then $n+1$ is a dominant
number of $T_t(R)$.}

\vspace{0.2cm}

{\bf Proof.} (1) If $R$ satisfies the $(l, n)^{op}$-condition, then
by Theorem 3.6, $\rfd_{T_t(R)} I^i(T_t(R))=\max\{\rfd_R I^i(R),
\rfd_R I^{i-1}(R)+1\}\leq l$ for any $0 \leq i\leq n-1$, which
implies that $T_t(R)$ satisfies the $(l+1, n)^{op}$-condition.
Conversely, by Theorem 3.6, we have that $\rfd_R I^i(R)\leq
\rfd_{T_t(R)} I^i(T_t(R))$ for any $i\geq 0$, so it is trivial that
$T_t(R)$ satisfies the $(l, n)^{op}$-condition implies so does $R$.

(2) If $n$ is a dominant number of $R$, then $\rfd_R I^i(R)< \rfd_R
I^n(R)$ for any $0 \leq i \leq n-1$. So by Theorem 3.6, for any $0
\leq i \leq n$, we have that $\rfd_{T_t(R)}
I^{n+1}(T_t(R))=\max\{\rfd_R I^{n+1}(R), \rfd_R I^{n}(R)+1\}\geq
\rfd_R I^{n}(R)+1>\max\{\rfd_R I^{i}(R), \rfd_R
I^{i-1}(R)+1\}=\rfd_{T_t(R)} I^{i}(T_t(R))$, which implies that
$n+1$ is a dominant number of $T_t(R)$. \hfill $\square$

\vspace{0.5cm}

{\bf Acknowledgements} This research was partially supported by the
Specialized Research Fund for the Doctoral Program of Higher
Education (Grant No. 20060284002), NSFC (Grant No. 10771095) and NSF
of Jiangsu Province of China (Grant No. BK2007517).

\end{document}